\documentclass[12pt]{article}
\usepackage{amsfonts}
\usepackage{mathrsfs}
\usepackage{amsmath,amssymb}
\baselineskip 5.5pt \lineskip 5.5pt \numberwithin{equation}{section}

\def\QED{\hfill$\Box$\par}

\def\SS{\mathcal {S}}
\def\HH{H}
\def\VV{\mathcal {V}ir}

\def\LL{\mathcal {L}}
\def\L{\mathcal {L}}

\def\cl{\centerline}

\def\rar{\longrightarrow}

\def\vs{\vspace*}

\def\ni{\noindent}

\def\C{\mathbb{C}}

\def\Z{\mathbb{Z}}
\def\adddot{$\!\!\!${\bf.}\ \ }

\newtheorem{theo}{Theorem}[section]

\newtheorem{lemm}[theo]{Lemma}

\begin{document}
\baselineskip 18pt

\cl{\bf{Leibniz Central Extension on the Twisted
Schr\"{o}dinger-Virasoro Algebra}\footnote {Supported by NSF grants
10471091, 10671027 of China.}}\vs{6pt}

\cl{Junbo Li$^{*,\dag)}$, Linsheng Zhu$^{\dag)}$} \cl{\small
$^{*)}$Department of Mathematics, Shanghai Jiao Tong University,
 Shanghai 200240, China}
\cl{\small $^{\dag)}$Department of Mathematics, Changshu Institute
of Technology, Changshu 215500, China} \cl{\small E-mail:
sd\_junbo@163.com, lszhu@cslg.edu.cn} \vs{6pt}

\noindent{\bf{Abstract.}} {In this paper we present all the Leibniz
2-cocycles of the twisted Schr\"{o}dinger-Virasoro algebra $\LL$},
which determine the second Leibniz cohomology group of $\LL$.
\noindent{\it{MR(2000) Subject Classification}: 17B05, 17B40, 17B65,
17B68.}\vs{20pt}

\cl{\bf1. \ Introduction}
\setcounter{section}{1}\setcounter{theo}{0}

Motivated by the research for the free Schr\"{o}dinger equations,
the original Schr\"{o}dinger-Virasoro Lie algebra was introduced in
\cite{H1}, in the context of non-equilibrium statistical physics.
The infinite-dimensional Lie algebra considered in this paper is
called the twisted Schr\"{o}dinger-Virasoro algebra, which is the
twisted deformation of the original Schr\"{o}dinger-Virasoro Lie
algebra. Both original and twisted sectors are closely related to
the Schr\"{o}dinger Lie algebras and the Virasoro Lie algebra, which
both play important roles in many areas of mathematics and physics
(e.g., statistical physics) and have been investigated in a series
of papers (see \cite{H2}, \cite{HU} and \cite{S1}--\cite{U}).

Now we give the definition of the Lie algebra $\LL$. A Lie algebra
$\LL$ is called a {\it twisted Schr\"{o}dinger-Virasoro Lie algebra}
(see \cite{RU}), if $\LL$ has the $\C$-basis
$$\{L_n,Y_n,M_n\,|\,n\in \Z\}$$
with the Lie brackets (others vanishing)
\begin{eqnarray}
&&[\,L_n,L_{n'}\,]=(n'-n)L_{n+n'},\ \ \
[\,L_n,M_p\,]=pM_{n+p},\label{LB1}\\[4pt]
&&[\,L_n,Y_m\,]=(m-\frac{n}{2})Y_{n+m},\ \ \
[\,Y_m,Y_{m'}\,]=(m'-m)M_{m+m'}.\label{LB2}
\end{eqnarray}
The twisted Schr\"{o}dinger-Virasoro Lie algebra has an
infinite-dimensional {\it twisted Schr\"{o}dinger subalgebra}
denoted by $\SS$ with the $\C$-basis $\{Y_n,M_n\,|\,n\in \Z\}$ and a
{\it Virasoro subalgebra} denoted by $\VV$ with the $\C$-basis
$\{L_n,C\,|\,n\in \Z\}$.

The Schr\"{o}dinger-Virasoro Lie algebras have recently drawn some
attentions in the literature. Particularly,  the sets of generators
provided by the cohomology classes of the cocycles for both original
and twisted sectors were presented in \cite{RU}, and the derivation
algebra and the automorphism group of the twisted sector were
determined in \cite{LS1}. Furthermore, vertex algebra
representations of these Lie algebras were constructed in \cite{U},
and the irreducible modules with finite-dimensional weight spaces
and indecomposable modules over them were considered in \cite{LS2}.

The main purpose of this paper is to determine the Leibniz
2-cocycles and further the second Leibniz cohomology group of
twisted Schr\"{o}dinger-Virasoro Lie algebra $\LL$ defined above. It
is well known that the 2-cocycles on Lie algebras play important
roles in the central extensions of Lie algebras, which can be used
to construct many infinite dimensional Lie algebras and in
particular, all 1-dimensional central extensions of $\LL$ are
determined by the 2-cohomology group of $\LL$. So do the Leibniz
2-cocycles and Leibniz cohomology groups of Lie algebras or further
Leibniz algebras. Therefore, there appeared a number of papers on
Leibniz 2-cocycles and Leibniz cohomology groups of infinite
dimensional Lie algebras and Leibniz algebras (see \cite{HPL},
\cite{L}--\cite{LP}, \cite{WT} and related references cited in
them). Now let's formulate our main results below.

We start with a brief definition. A {\it Leibniz algebra} $L$ over
$\C$ is a vector space equipped with a $\C$-bilinear map
$[\,\cdot\,,\,\cdot\,]: L\times L\rar L$ satisfying the Leibniz
identity:
\begin{eqnarray}\label{LJIdi}
[x,[y,z]]=[[x,y],z]-[[x,z],y],\ \ \ \forall\,\,x,y,z\in L.
\end{eqnarray}
It is easy to see that a Lie algebra must be a Leibniz algebra,
while a Leibniz algebra over $\C$ gives rise to be a Lie algebra if
$[x,x]=0$ for any $x\in L$.

Recall that a {\it Leibniz 2-cocycle} on $L$ is a $\C$-bilinear
function $\psi:L\times L\rar\C$ satisfying the Jacobian identity:
\begin{eqnarray}\label{2-cocy-JIdi}
\psi(x,[y,z])=\psi([x,y],z)-\psi([x,z],y)
\end{eqnarray}
for $x,y,z\in L$. In order to distinguish the Leibniz 2-cocycles
from the usual 2-cocycles (which are anti-symmetric in addition), we
call the usual ones Lie 2-cocycles. Denote by ${C_L^2}(L,\C)$ the
vector space of Leibniz 2-cocycles on $L$. For any $\C$-linear
function $f:L\rar\C$, one can define a Leibniz 2-cocycle $\psi_f$ as
follows
\begin{eqnarray}\label{cobLibco}
\psi_f(x,y)=f([x,y]),\ \ \forall\,\,x,y\in L,
\end{eqnarray}
which is called a {\it Leibniz 2-coboundary} or a {\it trivial
Leibniz 2-cocycle} on $L$. Denote by ${B_L^2}(L,\C)$ the vector
space of Leibniz 2-coboundaries on $L$. A Leibniz 2-cocycle $\phi$
is said to be {\it equivalent to} another Leibniz 2-cocycle $\psi$
if $\phi-\psi$ is trivial. For a 2-cocycle $\psi$, we denote by
$[\psi]$ the equivalent class of $\psi$. The quotient space
\begin{eqnarray}\label{dseco}
H_L^2(L,\C)\!=\!C_L^2(L,\C)/B_L^2(L,\C)\! =\{\mbox{the equivalent
classes of 2-cocycles}\},
\end{eqnarray}
is called the {\it second Leibniz cohomology group} of $L$.

\begin{theo}\adddot\label{mainth} The {\it second Leibniz cohomology group} of
$\L$, $\HH^2_L(\mathcal {L},\C)\cong\C$ is generated by the Virasoro
Leibniz cocycle (see (\ref{L2rede1})).
\end{theo}

Throughout the article, we denote by $\Z^*$ the set of all nonzero
integers and $\C^*$ the set of all nonzero complex numbers.

\vs{22pt}\cl{\bf2. Proof of the main
results}\setcounter{section}{2}\setcounter{theo}{0}\setcounter{equation}{0}

Let $\psi$ be any Leibniz 2-cocycle. Our main attempt or method is
to subtract all equivalent classes of the Leibniz 2-coboundaries on
$\LL$ from $\psi$. The proof of Theorem \ref{mainth} will be based
on several technical lemmas and divided into three cases with some
subcases.

According to the brackets (\ref{LB1})--(\ref{LB2}), it is easy to
see that the following identities hold:
\begin{eqnarray*}
L_n=\left\{\begin{array}{ll}
\frac{1}{n}[L_0,L_n] &{\rm if}\ \,n\neq0,\vs{6pt}\\
\frac{1}{2}[L_{-1},L_{1}] &{\rm if}\ \,n=0,
\end{array}\right.
M_n=\left\{\begin{array}{ll}
\frac{1}{n}[L_0,M_n] &{\rm if}\ \,n\neq0,\vs{6pt}\\
\,[L_{-1},M_{1}] &{\rm if}\ \,n=0,\end{array}\right.
Y_n=\left\{\begin{array}{ll}
\frac{1}{n}[L_0,Y_n] &{\rm if}\ \,n\neq0,\vs{6pt}\\
\frac{2}{3}[L_{-1},Y_{1}] &{\rm if}\ \,n=0.\end{array}\right.
\end{eqnarray*}
Define a $\C$-linear function $f:\LL\rightarrow\C$ as follows
\begin{eqnarray*}
&&f(Y_0)=\frac{2}{3}\psi(L_{-1},Y_{1}),\ \ \
f(L_0)=\frac{1}{2}\psi(L_{-1},L_{1}),\ \ \
f(M_0)=\psi(L_{-1},M_{1}),\ \ \ {\rm while}\\
&&f(Y_n)=\frac{1}{n}\psi(L_0,Y_n),\ \ \
f(L_n)=\frac{1}{n}\psi(L_0,L_n),\ \ \
f(M_n)=\frac{1}{n}\psi(L_0,M_n),\ \ \ {\rm if}\,\,n\in\Z^*.
\end{eqnarray*}
Let $\varphi=\psi-\psi_f$ where $\psi_f$ is defined in
(\ref{cobLibco}). One has
\begin{eqnarray}
&&\varphi(L_{-1},L_1)=\varphi(L_{-1},Y_1)=\varphi(L_{-1},M_1)=0,\label{reBa1}\\
&&\varphi(L_0,L_n)=\varphi(L_0,Y_n)=\varphi(L_0,M_n)=0,\ \
\,\,\forall\,\,n\in\Z^*.\label{reBa2}
\end{eqnarray}
\begin{lemm}\adddot\label{lemm1}
The following identities hold:
\begin{eqnarray}
&&\varphi(L_1,L_{-1})=\varphi(L_1,Y_{-1})=\varphi(L_1,M_{-1})=0,\label{rede1-2}\\
&&\varphi(L_n,L_0)=\varphi(Y_n,L_0)=\varphi(M_n,L_0)=0,\ \
\,\,\forall\,\,n\in\Z^*.\label{rede2}
\end{eqnarray}
\end{lemm}
\ni{\it Proof.} Replacing $\psi$ by $\varphi$ and $y$ by $x$ in
(\ref{2-cocy-JIdi}) simultaneously, one has
\begin{eqnarray}
\varphi(x,[x,z])+\varphi([x,z],x)=0,\
\,\forall\,\,x,z\in\LL.\label{LIcap1}
\end{eqnarray}
For any $n\in\Z^*$, replacing $(x,z)$ by
$(L_n,L_{-n}),\,\,(L_n,Y_{-n}),\,\,(L_n,M_{-n})$ in (\ref{LIcap1})
respectively, we obtain the following identity:
\begin{eqnarray*}
\varphi(L_{0},L_n)+\varphi(L_{n},L_0)
=\varphi(L_n,Y_0)+\varphi(Y_{0},L_n)=\varphi(L_{n},M_0)+\varphi(M_{0},L_n)=0,
\end{eqnarray*}
which together with (\ref{reBa2}) gives (\ref{rede2}).

Using the Jacobian identity on the three triples
$(L_1,L_{1},L_{-2}),\,(L_1,L_{1},Y_{-2}),\,(L_1,L_1,M_{-2})$
respectively, one has
\begin{eqnarray*}
\varphi(L_{-1},L_1)+\varphi(L_1,L_{-1})=
\varphi(L_{-1},Y_1)+\varphi(L_1,Y_{-1})=
\varphi(L_{-1},M_1)+\varphi(L_1,M_{-1})=0,
\end{eqnarray*}
which together with (\ref{reBa1}) gives (\ref{rede1-2}). Then the
lemma follows.\QED
\begin{lemm}\adddot\label{lemm2}
For any $m,n\in\Z$, one can write
\begin{eqnarray}
&&\varphi(L_n,L_{m})=\frac{n^3-n}{12}\delta_{m,-n}.\label{L2rede1}
\end{eqnarray}
\end{lemm}
\ni{\it Proof.} For any $m,n\in\Z$, applying the Jacobian identity
on the two triples $(L_{0},L_{-n},L_{n})$, $(L_{0},L_m,L_{n})$ and
$(L_{m},L_n,L_{-m-n})$, one immediately has
\begin{eqnarray*}
&&n(2\varphi(L_0,L_0)+\varphi(L_n,L_{-n})+\varphi(L_{-n},L_n)=0,\\
&&(m+n)\varphi(L_m,L_n)+(n-m)\varphi(L_{m+n},L_0)=0,
\end{eqnarray*}
which together with (\ref{rede2}), give
\begin{eqnarray}\label{cL1-1}
&&(m+n)\varphi(L_m,L_n)=\varphi(L_{0},L_0)=\varphi(L_n,L_{-n})+\varphi(L_{-n},L_n)=0.
\end{eqnarray}

For any $m,n\in\Z$, applying the Jacobian identity on the triple
$(L_m,L_n,L_{-m-n})$, we obtain
\begin{eqnarray}\label{corI-bd2}
&&(m-1)\varphi(L_{m+1},L_{-m-1}) =(m+2)\varphi(L_m,L_{-m}),
\end{eqnarray}
in which using induction on $n$, one can write
\begin{eqnarray}\label{corI-bd4}
&&\varphi(L_n,L_{-n})=\frac{n^3-n}{12},\ \ \ \forall\,\,n\in\C.
\end{eqnarray}
Then the lemma follows.\QED
\begin{lemm}\adddot\label{lemm3}
For any $m,n\in\Z$, one has
\begin{eqnarray*}
\varphi(M_m,M_{n})=\varphi(Y_m,M_{n})=0.
\end{eqnarray*}
\end{lemm}
\ni{\it Proof.} For any $m,n\in\Z$, applying the Jacobian identity
on the three triples $(L_{1},M_{-1},M_{0})$, $(Y_{0},Y_{m},M_{n})$
and $(Y_{-n},Y_{m},M_{n})$ respectively, one has
\begin{eqnarray*}
&&\varphi(M_{0},M_0)=m\varphi(M_{m},M_n)=(m+n)\varphi(M_{m-n},M_n)=0,
\end{eqnarray*}
which immediately gives
\begin{eqnarray*}
\varphi(M_{m},M_n)=0\ \ \ \forall\,\,m,n\in\Z.
\end{eqnarray*}

For any $m,n\in\Z$, applying the Jacobian identity on the three
triples $(L_{1},Y_{-1},M_{0})$, $(L_{0},Y_{n},M_{-n})$ and
$(L_{m},Y_{0},M_{n})$, one has
\begin{eqnarray}
&&n\big(\varphi(M_{-n},Y_{n})+\varphi(Y_{n},M_{-n})\big)=0,\label{L22-10}\\
&&2n\varphi(M_{m+n},Y_0)+m\varphi(Y_{m},M_n)=0,\label{L22-11}
\end{eqnarray}
which immediately gives
\begin{eqnarray}\label{L22-1}
n\varphi(M_{n},Y_0)=\varphi(Y_{n},M_0)=0\ \ \,\forall\,\,n\in\Z.
\end{eqnarray}
and further implies
\begin{eqnarray}\label{L22-2}
m(m+n)\varphi(Y_{m},M_n)=0\ \ \,\forall\,\,n\in\Z.
\end{eqnarray}

For any $n\in\Z$, applying the Jacobian identity on the two triples
$(L_{n},Y_{0},M_{-n})$ and $(Y_{n},Y_{0},Y_{-n})$, one has
\begin{eqnarray*}
&&n\big(\varphi(Y_{n},M_{-n})-2\varphi(M_{0},Y_{0})\big)=0,\\
&&n\big(\varphi(M_{n},Y_{-n})-\varphi(Y_{n},M_{-n})
-2\varphi(M_{0},Y_{0})\big)=0,
\end{eqnarray*}
which gives
\begin{eqnarray*}
\varphi(M_{n},Y_{-n})=2\varphi(Y_{n},M_{-n})=4\varphi(M_{0},Y_{0})\
\ \,\forall\,\,n\in\Z^*,
\end{eqnarray*}
and together with (\ref{L22-10}), further infers
\begin{eqnarray*}
\varphi(M_{n},Y_{-n})=\varphi(Y_{n},M_{-n})=0\ \
\,\forall\,\,n\in\Z.
\end{eqnarray*}
Then the lemma follows.\QED
\begin{lemm}\adddot\label{lemm4}
For any $m,n\in\Z$, one has
\begin{eqnarray*}
&&\varphi(L_n,M_{m})=\varphi(M_{n},L_{m})=nc_1\delta_{m,-n},\ \
\varphi(Y_n,Y_{m})=2nc_1\delta_{m,-n}
\end{eqnarray*}
for some constant $c_1\in\C$.
\end{lemm}
\ni{\it Proof.} For any $m,n\in\Z$, applying the Jacobian identity
on the three triples $(L_m,L_0,M_n)$, $(L_{1},L_{-1},M_{0})$ and
$(L_m,L_0,M_n)$ respectively, one has
\begin{eqnarray*}
&&n\varphi(L_{0},M_{m+n})+n\varphi(M_n,L_m)-m\varphi(L_{m},M_n)=0,\\
&&\varphi(L_{0},M_0)=(m+n)\varphi(L_m,M_n)+n\varphi(M_{m+n},L_0)=0,
\end{eqnarray*}
which together with (\ref{reBa2}) and (\ref{rede2}), immediately
give ($\forall\,\,m,n\in\Z$)
\begin{eqnarray}
&&n\varphi(M_n,L_m)=m\varphi(L_{m},M_n),\label{lem4re0-1}\\
&&(m+n)\varphi(L_m,M_n)=\varphi(M_{0},L_0)
=\varphi(L_{0},M_0)=0.\label{lem4re0-2}
\end{eqnarray}
For any $n\in\Z$, applying the Jacobian identity on the triple
$(L_n,L_1,M_{-n-1})$, one has
\begin{eqnarray*}
(n-1)\varphi(L_{n+1},M_{-n-1})
=(n+1)\big(\varphi(L_n,M_{-n})-\varphi(L_1,M_{-1})\big),
\end{eqnarray*}
which gives
\begin{eqnarray}\label{lem4re1}
\varphi(L_n,M_{-n})=\left\{\begin{array}{ll}
\frac{1}{2}n(n-1)\varphi(L_2,M_{-2})+n(2-n)\varphi(L_1,M_{-1})
&{\rm if}\ \,n\geq0,\vs{6pt}\\
\frac{1}{6}n(n-1)\varphi(L_{-2},M_{2})+\frac{1}{3}n(n+2)\varphi(L_1,M_{-1})
&{\rm if}\ \,n\leq-1.
\end{array}\right.
\end{eqnarray}
Applying the Jacobian identity on the triple $(L_1,L_{-2},M_{1})$,
one has
\begin{eqnarray*}
3\varphi(L_{-1},M_{1})=\varphi(L_{-2},M_{2})-\varphi(L_1,M_{-1}),
\end{eqnarray*}
which together with the case $n=-1$ of (\ref{lem4re1}) gives
\begin{eqnarray*}
\varphi(L_{-2},M_{2})=3\varphi(L_{2},M_{-2})-8\varphi(L_1,M_{-1}),
\end{eqnarray*}
and together with which, (\ref{lem4re1}) becomes
\begin{eqnarray}\label{lem4re2}
\varphi(L_n,M_{-n})=n(n-1)c_1+n(n-2)c_2\ \ \,\forall\,\,n\in\Z,
\end{eqnarray}
by denoting $\varphi(L_1,M_{-1})=-c_2,\,\varphi(L_2,M_{-2})=2c_1$.
Using (\ref{lem4re0-1}), (\ref{lem4re0-2}) and (\ref{lem4re2}), we
obtain
\begin{eqnarray}\label{lem4re22}
\varphi(M_{n},L_{-n})=-n(n+1)c_1-n(n+2)c_2\ \ \,\forall\,\,n\in\Z,
\end{eqnarray}

For any $m,n\in\Z$, applying the Jacobian identity on the two
triples $(L_{0},Y_{m},Y_{n})$ and $(L_{m},Y_{0},Y_{n})$, one has
\begin{eqnarray}
&&(n-m)\varphi(L_0,M_{m+n})+n\varphi(Y_n,Y_m)-m\varphi(Y_m,Y_n)=0,\nonumber\\
&&2n\varphi(L_m,M_n)+m\varphi(Y_m,Y_n)+(2n-m)\varphi(Y_{m+n},Y_0)=0,\label{lem5re0-0}
\end{eqnarray}
which together with (\ref{reBa2}) and (\ref{lem4re0-2}), immediately
give
\begin{eqnarray}
&&n\varphi(Y_n,Y_m)=m\varphi(Y_m,Y_n)\ \
\,\forall\,\,m,n\in\Z,\label{lem5re0-2}
\end{eqnarray}
which further infers
\begin{eqnarray}
&&n\varphi(Y_n,Y_0)=n\big(\varphi(Y_n,Y_{-n})+\varphi(Y_{-n},Y_n)\big)=0\
\ \,\forall\,\,n\in\Z.\label{lem5re1-0}
\end{eqnarray}
Then recalling (\ref{lem4re2}), (\ref{lem4re22}), (\ref{lem5re0-2})
and (\ref{lem5re1-0}), one can rewrite (\ref{lem5re0-0}) as
\begin{eqnarray*}
2n\varphi(L_m,M_n)+n\varphi(Y_n,Y_m)+3n\varphi(Y_{m+n},Y_0)=0,\label{lem5re2-1}
\end{eqnarray*}
from which one can obtain the following identity:
\begin{eqnarray}\label{laLM4a2}
\varphi(Y_n,Y_{m})=
-2n\big((n+1)c_1+(n+2)c_2\big)\delta_{m,-n}-3\varphi(Y_{0},Y_0)\delta_{m,-n}\,\
\ \forall\,\,m,n\in\Z^*.
\end{eqnarray}
Then noticing the following identities:
\begin{eqnarray*}
&&\varphi(Y_n,Y_{-n})=
-2n\big((n+1)c_1+(n+2)c_2\big)-3\varphi(Y_{0},Y_0)\ \ \ n\in\Z^*,\\
&&\varphi(Y_{-n},Y_{n})=
2n\big((-n+1)c_1+(-n+2)c_2\big)-3\varphi(Y_{0},Y_0)\ \ \ n\in\Z^*,
\end{eqnarray*}
and recalling (\ref{lem5re1-0}), one has
\begin{eqnarray*}
&&2n\big((-n+1)c_1+(-n+2)c_2\big)-3\varphi(Y_{0},Y_0)\\
&&=2n\big((n+1)c_1+(n+2)c_2\big)+3\varphi(Y_{0},Y_0),
\end{eqnarray*}
which gives
\begin{eqnarray*}
3\varphi(Y_{0},Y_0)+2(c_1+c_2)n^2=0\ \ \ \forall\,\,n\in\Z^*,
\end{eqnarray*}
and further forces
\begin{eqnarray}\label{laLM4a21}
\varphi(Y_{0},Y_0)=c_1+c_2=0.
\end{eqnarray}
Then using (\ref{laLM4a21}), one can rewrite (\ref{lem4re2}),
(\ref{lem4re22}) and (\ref{laLM4a2}) respectively as follows:
\begin{eqnarray*}
&&\varphi(Y_n,Y_{-n})=2\varphi(L_n,M_{-n})=2\varphi(M_{n},L_{-n})=2nc_1\
\ \,\forall\,\,n\in\Z.
\end{eqnarray*}
Then the lemma follows.\QED

\begin{lemm}\adddot\label{lemm6}
For any $m,n\in\Z$, one has
\begin{eqnarray*}
\varphi(L_m,Y_{n})=\varphi(Y_{n},L_m)=0.
\end{eqnarray*}
\end{lemm}
{\it Proof.} For any $m,n\in\Z$, applying the Jacobian identity on
the three triples $(L_m,L_0,Y_0)$, $(L_m,L_0,Y_n)$, $(L_m,L_n,Y_0)$
and $(L_0,L_m,Y_n)$, one has
\begin{eqnarray*}
&&m\big(2\varphi(L_m,Y_0)-\varphi(Y_m,L_0)\big)=0,\\
&&2(m+n)\varphi(L_m,Y_n)-(m-2n)\varphi(Y_{m+n},L_0)=0,\\
&&n\varphi(L_m,Y_n)+2(n-m)\varphi(L_{m+n},Y_0)+m\varphi(Y_m,L_n)=0,\\
&&2m\varphi(L_m,Y_n)-2n\varphi(Y_n,L_m)+(m-2n)\varphi(L_0,Y_{m+n})=0,
\end{eqnarray*}
which further can be rewritten as follows \big(recalling
(\ref{reBa2}) and (\ref{rede2})\big):
\begin{eqnarray}
&&\varphi(L_m,Y_0)=0\ \ \forall\,\,m\in\Z^*,\ \
\varphi(L_m,Y_n)=0\ \ \forall\,\,m+n\in\Z^*,\label{LM5LY0}\\
&&\varphi(L_{-n},Y_n)=\varphi(Y_{-n},L_n)
-3\varphi(L_{0},Y_0)\ \ \forall\,\,n\in\Z^*,\label{LM5LY1}\\
&&\varphi(L_{-n},Y_n)=-\varphi(Y_n,L_{-n})-\frac{3}{2}\varphi(L_0,Y_{0})\
\ \forall\,\,n\in\Z^*.\label{LM5LY2-1}
\end{eqnarray}
The identities (\ref{LM5LY1}), (\ref{LM5LY2-1}) together with
(\ref{reBa1}) and (\ref{rede1-2}) force
\begin{eqnarray}
&&\varphi(Y_{-1},L_1)=3\varphi(L_{0},Y_0)=0,\label{LM5LY1-1}
\end{eqnarray}
which together with (\ref{LM5LY1}) and (\ref{LM5LY2-1}) gives
\begin{eqnarray}
\varphi(L_{n},Y_{-n})+\varphi(L_{-n},Y_n)=0.\label{LM5LY88}
\end{eqnarray}

For any $n\in\Z$, applying the Jacobian identity on the triple
$(L_n,L_1,Y_{-n-1})$ and recalling (\ref{LM5LY1-1}), one has
\begin{eqnarray*}
2(n-1)\varphi(L_{n+1},Y_{-1-n})=(2n+3)\varphi(L_n,Y_{-n}),
\end{eqnarray*}
which gives (\,by using induction on $n$)
\begin{eqnarray}\label{LM5LY3-1}
\varphi(L_n,Y_{-n})=0\ \ \ \forall\,\,n\in\Z.
\end{eqnarray}
which together with (\ref{LM5LY1}) further forces
\begin{eqnarray}\label{LM5LY3-6}
\varphi(Y_n,L_{-n})=0\ \ \ \forall\,\,n\in\Z.
\end{eqnarray}
Then this lemma follows from (\ref{LM5LY0}), (\ref{LM5LY3-1}) and
(\ref{LM5LY3-6}).\QED
\begin{lemm}\adddot\label{lemm66}
One can suppose
\begin{eqnarray*}
c_1=0.
\end{eqnarray*}
\end{lemm}
{\it Proof.} Define another $\C$-linear function
$g:\LL\rightarrow\C$ as follows
\begin{eqnarray*}
&&g(M_0)=c_1,\ \ {\mbox{other components vanishing}}.
\end{eqnarray*}
Still denote $\varphi-\psi_g$ by $\varphi$ where $\psi_g$ is defined
in (\ref{cobLibco}). One has
\begin{eqnarray}
c_1=\varphi(L_{-1},M_1)=0.\label{reBa18}
\end{eqnarray}

\ni{\it Proof of Theorem \ref{mainth}}\ \,The theorem follows by the
series of lemmas from the second one to the last one.

\end{document}